\documentclass[12pt]{article}

\usepackage{a4wide}
\usepackage{amsmath,amssymb,amsfonts,amsthm}

\begin{document}

\title{Proof of Two Conjectures on Correlation Inequalities\\ for One Class of Monotone
Functions}

\author{V. Blinovsky\footnote{The author was supported by FAPESP (2012/13341-8, 2013/07699-0) and
  NUMEC/USP (Project MaCLinC/USP)}}

\date{\small
 Instituto de
Matem\'atica e Estat\'{\i}stica, Universidade de S\~ao Paulo, Brazil\\
Institute for Information Transmission Problems, Moscow, Russia\\ {\tt
vblinovs@yandex.ru}\\[2ex] February 10, 2013}

\maketitle

\begin{abstract}
We prove two correlation inequalities conjectured in~\cite{1} for functions that are
linear combinations of unimodal Boolean monotone nondecreasing functions.
\end{abstract}

\bigskip
Assume that $f_1,\ldots,f_n$ are nonnegative nondecreasing functions
$2^X\to\mathbb{R}$, where $X=[m]$ for some finite $m\ge 1$. The expectation of a
random variable $f\colon 2^X\to\mathbb{R}$ with respect to $\mu$ is denoted by
$\langle f\rangle_\mu$. For a subset $\delta\subset [n]$, define
$$
E_\delta =\biggl\langle\,\prod_{i\in\delta} f_i\biggr\rangle_{\!\mu}.
$$
Let
$$
\sigma =\{\sigma_1,\ldots,\sigma_\ell\}
$$
be a partition of $[n]$ into disjoint subsets. Define
$$
E_\sigma =\prod_{i=1}^{\ell} E_{\sigma_i}.
$$
Let $\lambda_i =|\sigma_i|$. We have $\smash[t]{\sum\limits_{i=1}^{\ell} \lambda_i}
=n$. Let $\lambda (\sigma)=(\lambda_1,\ldots,\lambda_{\ell})$ and $\lambda_1\ge
\ldots \ge\lambda_{\ell}$. For a partition~$\lambda$ of~$n$, $n=\lambda_1
+\ldots+\lambda_{\ell}$, define
\begin{equation}\label{ee5}
E_\lambda =\sum_{\sigma :\: \lambda (\sigma)=\lambda}E_{\sigma}.
\end{equation}
A unimodal monotone nondecreasing Boolean function $f$ on $2^X$ is a function that
takes values in $\{0,1\}$ and has the properties that $f(A)\ge f(B)$ for $B\subset A$
and there is a unique $C\in 2^X $ such that $f(C)=1$ and $f(C_1)=0$ for an arbitrary
$C_1\subset C$.

In~\cite{1} it was conjectured that, if the probability measure $\mu$ on $2^X$
satisfies the FKG condition
\begin{equation}\label{e1}
\mu (A\wedge B)\mu (A\vee B)\ge \mu (A)\mu (B),
\end{equation}
then
\begin{equation}\label{ee3}
E_{n,\mu} (f_1,\ldots,f_n)=\sum_{\lambda\vdash n}c_{\lambda}E_{\lambda} \ge 0,
\end{equation}
where
$$
c_{\lambda}=(-1)^{\ell +1}\prod_{i=1}^{\ell}(\lambda_i -1)!.
$$
In~\cite{1} this conjecture was proved for a product measure $\mu$ and for functions
that are linear combinations of unimodal monotone nondecreasing Boolean functions
with nonnegative coefficients (actually, in~\cite{1}, under these assumptions
conjecture~\eqref{e2} was proved and it was shown that conjecture~\eqref{e2} implies
inequalities~\eqref{ee3}).

Now we prove inequalities~\eqref{ee3} without assuming that $\mu$ is a product
measure, i.e., for an arbitrary probability measure that satisfies the FKG condition,
and show that inequalities~\eqref{ee3} imply the validity of conjecture~\eqref{e2}.

Since the functional $E_n (\cdot)$ is linear in each function $f_i$, it is sufficient
to prove~\eqref{ee3} for monotone nondecreasing Boolean unimodal functions. Let $f_i$
have support $A_i$.

\smallskip
\textbf{Lemma.}
\emph{If\/ $A,B,C\subset 2^{X}$ are supports of unimodal monotone nondecreasing Boolean
functions\/ \textup(upsets with one minimal element\textup) and\/ $\mu$ satisfies the FKG
condition\textup, then the conditional probability\/ $\mu_C$ satisfies the inequality}
$$
\mu_C (A\cap B)\ge \mu_C (A)\mu_C (B).
$$

\smallskip
The lemma easily follows from the fact that the measure $\mu$ satisfies the
inequality
$$
\mu (C)\mu (A\cap B\cap C)\ge \mu (A\cap C)\mu (B\cap C).
$$
This in turn is a consequence of the four function theorem~\cite{2}.

It is convenient to use the following notation:
\begin{gather*}
E_{n,\mu} (f_1,\ldots,f_n)=E_{n,\mu} (A_1,\ldots,A_n)=E_{n,\mu} (\{A_i,\: i\in
[n]\}),\\ \langle f_1 f_2 \ldots f_n \rangle_\mu =\mu \biggl(\,\bigcap_{i=1}^n A_i
\biggr) =\mu (\{A_i,\: i\in [n]\}).
\end{gather*}

Consider the decomposition
\begin{equation}\label{er6}
E_{n,\mu} (\{A_i,\: i\in [n]\})=(n-1)!\,\mu (\{A_i,\: i\in [n]\}) -I_{n,\mu}
(\{A_i,\: i\in [n]\}),
\end{equation}
where the first term comes from the sum in formula~\eqref{ee3} for a partition
$\lambda$ with $\ell =1$ and $I_{n,\mu}$ is a correction term added to make the
formula valid. We are going to prove inequality~\eqref{ee3} for Boolean unimodal
functions in the general way. Let $\sigma (A_n)\in\sigma$, $n\in\sigma (A_n)$, and
$k\le n-1$. Define
\begin{equation}\label{er7}
\begin{gathered}
E^k_{n,\mu} (\{A_i,\: i\in [n]\}) = \sum_{\lambda\vdash n}c_\lambda
E_\lambda^k,\qquad E^k_\lambda =\sum_{\sigma:\:\lambda
(\sigma)=\lambda}E^k_\sigma,\\[-2pt]
\begin{aligned}
E_\sigma^k &=\mu^{|[k]\cap\sigma (A_n)|}(A_n)\mu (\{A_i,\: i\in\sigma (A_n)\})\\
&\quad\strut\times\prod_{\sigma_i\in\sigma\setminus \sigma (A_n)}\mu^{||\sigma_i\cap
[k]|-1|_+} (A_n) \mu (\{A_i\cap A_n,\: i\in\sigma_ i\cap [k]\} \{A_i,\: i\in\sigma_i
\setminus [k]\}),
\end{aligned}
\end{gathered}
\end{equation}
where
$$
|x|_+=
\begin{cases}
0,& x\le0,\\ x,& x>0.
\end{cases}
$$
It is clear that
$$
E_{n,\mu} (\{A_i,\: i\in [n]\})=E^0_{n,\mu} (\{A_i,\: i\in [n]\}).
$$
It is important that, in the notation $\mu (\{A_i,\: i\in Y\})$, each $A_i$ has its
own transformations, but each new $A_i$ still keeps its place, even though the
product $\bigcap A_i$ can be written in terms of $A_i$ in different ways. Also we
assume that all $A_i$ in the notation $\mu (\{A_i,\: i\in Y\})$ are ordered in a
natural way.

In the Appendix we write most of the formulas from the main part of the text in a
particular case of $n=3$ to make them more transparent. The case of an arbitrary $n$
is similar to the case of $n=3$.

Next we prove that
\begin{equation}\label{er8}
E^k_{n,\mu} (\{A_i,\: i\in [n]\})\ge 0.
\end{equation}
We will use induction assuming that inequality~\eqref{er8} is true for $n-1$ and for
all $k\le n-2$. It is obviously true for $n=2$.

Actually, first we are going to prove that
\begin{multline}\label{er9}
E^{k}_{n,\mu} (\{A_i,\: i\in [n]\})\ge E^{k+1}_n (\{A_i,\: i\in [n]\})\ge
E^{n-1}_{n,\mu} (\{A_i,\: i\in [n]\})\\= \mu^{n}(A_n)E_{n,\mu_{A_n}}(\{A_i,\: i\in
[n]\}).\quad
\end{multline}
Equality in~\eqref{er9} is clear from the definition of $E^k_{n,\mu}$.

Indeed,
$$
\begin{aligned}
E_{n,\mu}^{n-1}(\{A_i,\: i\in [n]\})&=\mu^n (A_n)\mu (\{A_i,\: i\in\sigma (A_n)\})
\prod_{\sigma_i\in\sigma\setminus \sigma (A_n)}\mu^{-\ell} (A_n) \mu (\{A_i\cap
A_n,\: i\in\sigma_ i\})\\ &=\mu^n (A_n) \frac{\mu (\{A_i,\: i\in\sigma (A_n)\})}{\mu
(A_n)}\prod_{\sigma_i\in\sigma\setminus \sigma (A_n)}\frac{\mu (\{A_i\cap A_n,\:
i\in\sigma_ i\})}{\mu (A_n)}\\ &=\mu^{n}(A_n)E_{n,\mu_{A_n}}(\{A_i,\: i\in [n]\}).
\end{aligned}
$$

Consider the decomposition
$$
E^k_{n,\mu}(\{A_i,\: i\in [n]\})=\mu^{k}(A_n)\mu (\{A_i,\: i\in [n]\})- I_{n,\mu}^k
(\{A_i,\: i\in [n]\}).
$$
We need to show the validity of the following inequality:
\begin{equation}\label{er0}
I_{n,\mu}^{k+1} (\{A_i,\: i\in [n]\})\ge \mu (A_n)I_{n,\mu}^k (\{A_i,\: i\in [n]\}).
\end{equation}
We have
$$
\begin{aligned}
\mu (A_n)I_{n,\mu}^k (\{A_i,\: i\in [n]\}) &=\mu (A_n)\sum_{Z\subset [k+1,n-1],
\{k+1\}\in Z}(|Z|-1)!\,\mu\left(\{A_i,\: i\in Z\} \right)\\ &\times E_{n-|Z|}^k
(\{A_i\cap A_n,\: i\in [k]\}, \{A_i,\: i\in [k+1,n]\setminus Z\}) +\mu (A_n)B\\
&\le\sum_{Z\subset [k+1,n-1], \{k+1\}\in Z}(|Z|-1)!\, \mu\left(\{A_i,\: i\in
Z\setminus\{k+1\}\}, A_{k+1}\cap A_n \right)\\ &\times E_{n-|Z|}^k (\{A_i\cap A_n,\:
i\in [k]\}, \{A_i,\: i\in [k+1,n]\setminus Z\}) +\mu (A_n)B\\ &=I_{n,\mu}^{k+1}
(\{A_i,\: i\in [n]\}).
\end{aligned}
$$
Here $B$ is the correction term which makes the first equality valid. To obtain the
inequality in the last chain of relations, we use the FKG inequality and the
induction assumption. Thus inequality~\eqref{er0} follows. Next we repeat all the
procedure with
$$
E_{n,\mu_{A_n}}(\{A_i,\: i\in [n]\})
$$
and new $\smash{A^\prime_i} = A_i\cap A_n$ using the conditional measure
$\smash[b]{\mu_{A^\prime_{n-1}}}=\smash[b]{\mu_{A_{n-1}\cap A_n}}$. Repeating the
same procedure $n$ times, we come to the inequality ($A=\bigcap\limits_{i=1}^n A_i$)
$$
E^{k}_{n,\mu} (\{A_i,\: i\in [n]\})\ge \mu^n (A)E_n (1,\ldots,1)=0.
$$
This proves inequality~\eqref{er8} (and also \eqref{ee3}).

Next we consider the set of formal series $P[[t]]$ whose coefficients are monotone
nondecreasing nonnegative functions on $2^X$. Then $p(A)=p_1 (A)t +p_2 (A)t^2 +\ldots
\in P[[t]]$.

In~\cite{1} the following conjecture was also formulated.

\smallskip
\textbf{Conjecture.}
\emph{For the FKG probability measure\/ $\mu$\textup, the following inequality is true}:
\begin{equation}\label{e2}
1-\prod_{A\in 2^X}(1-p(A))^{\mu (A)}\ge 0.
\end{equation}

\smallskip
Inequality~\eqref{e2} is understood as nonnegativeness of coefficients of the formal
series obtained by series expansion of the product on the left-hand side of this
inequality.

We will prove that inequality~\eqref{e2} follows from
\begin{equation}\label{e56}
E_n (f_1,\ldots, f_n)\ge 0
\end{equation}
for all $n$, and hence it is sufficient to prove the last inequality and then
inequality~\eqref{e2} follows.

Let us make some transformations in the expression on the left-hand side
of~\eqref{e2}. We have
\begin{align*}
1-\prod_{A\in 2^X}(1-p(A))^{\mu (A)}&=1-\exp\left\{\left\langle \ln
(1-p)\right\rangle_\mu\right\}\\[-5pt]
&=1-\exp\left\{-\sum_{i=1}^{\infty}\frac{1}{i}\langle p^i \rangle_{\mu}\right\}\\ &=
\sum_{j=1}^{\infty}\frac{(-1)^{j+1}}{j!}\left(\sum_{i=1}^{\infty}\frac{1}{i} \langle
p^i \rangle_{\mu}\right)^j\\[-10pt] &=
\sum_{j=1}^{\infty}\frac{(-1)^{j+1}}{j!}\sum_{\{q_s\} :\:\sum q_s
=j}\binom{j}{q_1,\ldots, q_j}\sum_{\{i_s\}} \frac {\prod\limits_{s=1}^j \langle
p^{i_s} \rangle^{q_s}_\mu}{(i_1)^{q_1} (i_2)^{q_2}\ldots (i_j)^{q_j}}.
\end{align*}
Next recall that the number of partitions of $n$ with a given set $\{q_i\}$ of
occurrences of $i$ is
$$
\frac{n!}{\prod\limits_i (i!)^{q_i}q_i!}.
$$
Continuing the last chain of identities and using the last formula, we obtain
\begin{equation}\label{e4}
\begin{aligned}[b]
1-\prod_{A\in 2^X}(1-p(A))^{\mu (A)} &=\sum_{n=1}^{\infty}\frac{1}{n!}
\sum_{\lambda\vdash n}\,\sum_{\sigma:\: \lambda (\sigma)=\lambda} (-1)^{\sum\limits_i
q_i +1}\,\frac{n!\,\prod ((i-1)!)^{q_i}}{\prod\limits_i (i!)^{q_i} q_i!}
\prod\limits_i \langle p^{i} \rangle^{q_i}_\mu\\ &=\sum_{n=1}^{\infty}\frac{1}{n!}
\sum_{\lambda\vdash n} (-1)^{\ell (\lambda)+1}\prod_i (\lambda_i -1)! \sum_{\sigma:\:
\lambda (\sigma)=\lambda} E_\sigma (p,\ldots,p)\\ &= \sum_{n=1}^{\infty}\frac{1}{n!}
\sum_{\lambda\vdash n} c_\lambda E_\lambda (p, \ldots,p)\\
&=\sum_{n=1}^{\infty}\frac{1}{n!}E _n (p,\ldots,p).
\end{aligned}
\end{equation}
Hence, to prove conjecture~\eqref{e2}, now we need to show that
\begin{equation}\label{et}
E_n (p,\ldots,p)\ge 0.
\end{equation}
But the coefficients of the formal series $E_n (p,\ldots,p)$ are sums of $E_n
(p_{i_1},\ldots, p_{i_n})$ for multisets $\{i_1,\ldots,i_n\}$. This completes the
proof that inequality~\eqref{e2} follows from inequality~\eqref{e56} under the same
conditions on $\mu$.

Therefore, we also prove inequality~\eqref{e2}.

\appendix

Here, by example of the case of $n=3$, we demonstrate some relations that we used and
proved above. From the definitions we obtain the equalities
\begin{align*}
E_{3,\mu}(A_1,A_2,A_3)&=2\mu (A_1\cap A_2\cap A_3) -\Bigl(\mu (A_1)\mu(A_2\cap A_3)\\
&\quad\strut+\mu (A_2)\mu(A_1\cap A_3) +\mu (A_3)\mu(A_1\cap A_2) -\mu (A_1)\mu
(A_2)\mu (A_3)\Bigr),\\ E^1_{(123)}=E^1_{(12)(3)} &=\mu (A_3)\mu(A_1\cap A_2\cap
A_3),\\ E^1_{(1)(23)} &=\mu (A_2\cap A_3) \mu(A_1\cap A_3),\\
E^1_{(13)(2)}=E_{(1)(2)(3)}^1 &=\mu (A_3)\mu (A_2)\mu(A_1\cap A_3).
\end{align*}
Thus,
$$
E^1_{3,\mu} (A_1,A_2,A_3)=\mu (A_3)\mu(A_1\cap A_2\cap A_3) -\mu(A_2\cap A_3)
\mu(A_1\cap A_3).
$$
Next,
\begin{gather*}
E^2_{(123)}=E^2_{(12)(3)} =\mu^2 (A_3)\mu(A_1\cap A_2\cap A_3),\\
E^2_{(1)(23)}=E^2_{(13)(2)} =E^2_{(1)(2)(3)}=\mu (A_3)\mu (A_2 \cap A_3) \mu(A_1\cap
A_3).
\end{gather*}
Thus,
$$
E^2_{3,\mu}(A_1,A_2,A_3)=\mu^2 (A_3)\mu(A_1\cap A_2\cap A_3) -\mu (A_3)\mu(A_2\cap
A_3) \mu(A_1\cap A_3)
$$
and
$$
I^2_{3,\mu} (A_1,A_2,A_3)=\mu (A_3) I^1_{3,\mu} (A_1,A_2,A_3).
$$
Validity of the inequality
$$
I^1_{3,\mu}(A_1,A_2,A_3)\ge\mu (A_3) I^0_{3,\mu} (A_1,A_2,A_3)
$$
follows from the chain of relations
\begin{align*}
\lefteqn{\mu (A_3) \Bigl( \mu ( A_1)\mu( A_2\cap A_3) +\mu (A_2)\mu(A_1\cap A_3) +\mu
(A_3)\mu(A_1\cap A_2) -\mu (A_1)\mu (A_2)\mu (A_3)\Bigr)}\\ &=\mu (A_3)\mu
(A_1)\Bigl(\mu(A_2\cap A_3) -\mu (A_2)\mu (A_3)\Bigr) +\mu (A_3)(\mu (A_2)\mu(A_1\cap
A_3) +\mu (A_3)\mu(A_1\cap A_2))\\ &\le \mu(A_1\cap A_3 )\Bigl(\mu(A_2\cap A_3) -\mu
(A_2)\mu (A_3)\Bigr)+\mu (A_3)(\mu (A_2)\mu(A_1\cap A_3) +\mu(A_1\cap A_2\cap A_3))\\
&=\mu(A_1\cap A_3) \mu(A_2\cap A_3) +\mu (A_3) \mu(A_1\cap A_2 \cap A_3)
=I^1_{3,\mu}(A_1,A_2,A_3).
\end{align*}
Here, to prove the inequality in the last chain of relations, we use the FKG
inequality.

\medskip
The work was done during the author's visit to the Institute of Mathematics and
Statistics of the University of S\~ao Paulo. The author is grateful to the Institute
for hospitality.

\end{document}